\documentclass[12pt]{article}
\usepackage{amssymb,amsmath,graphicx,amsthm,mathrsfs,bm}
\usepackage{color}
\newtheorem{Theorem}{Theorem}[section]
\newtheorem{Lemma}[Theorem]{Lemma}

\numberwithin{equation}{section}

\newcommand{\lc}
{\mathrel{\raise2pt\hbox{${\mathop<\limits_{\raise1pt\hbox
{\mbox{$\sim$}}}}$}}}

\newcommand{\gc}
{\mathrel{\raise2pt\hbox{${\mathop>\limits_{\raise1pt\hbox{\mbox{$\sim$}}}}$}}}

\newcommand{\ec}
{\mathrel{\raise2pt\hbox{${\mathop=\limits_{\raise1pt\hbox{\mbox{$\sim$}}}}$}}}

\begin{document}

\title{  A Game Problem for Heat Equation }

\author{Lijuan Wang\thanks{School of
Mathematics and Statistics, Computational Science Hubei Key
Laboratory, Wuhan University, Wuhan 430072, China; e-mail:
ljwang.math@whu.edu.cn. The first author was supported by the
National Natural Science Foundation of China under grant 11771344.},
\quad Donghui Yang\thanks{School of Mathematics and Statistics,
School of Information Science and Engineering, Central South
University, Changsha 410075, China; e-mail: donghyang@139.com. The
second author was supported by  National Natural Science Foundation
of China, China Postdoctoral Science Foundation under grant
2017M610503 and Central South University Postdoctoral Science
Foundation,} \quad and \quad Zhiyong Yu\thanks{Corresponding author.
School of Mathematics, Shandong university, Jinan 250100, China;
email: yuzhiyong@sdu.edu.cn. The third author was supported by the
National Natural Science Foundation of China under grant 11471192
and the Fundamental Research Funds of Shandong University under
grant 2017JC016.}}

\date{}
\maketitle
\begin{abstract}
In this paper, we consider a two-person game problem governed by a linear heat equation.
The existence of Nash equilibrium for this problem is considered. Moreover, the bang-bang
property of Nash equilibrium is discussed.
\end{abstract}

\medskip


\medskip


\section{Introduction}

In this paper, we assume that $\Omega$ is a bounded domain in $\mathbb{R}^d, d\geq 1$, with
a $C^2$ boundary $\partial\Omega$. Let $T$ be a given positive constant,
$\omega_1$ and $\omega_2$ be two open and nonempty subsets
of $\Omega$ with $\omega_1\cap \omega_2=\emptyset$. Let $\chi_{\omega_1}$ and $\chi_{\omega_2}$ be two
characteristic functions of the sets $\omega_1$ and $\omega_2$, respectively. The controlled linear heat
equation under consideration is as follows:
\begin{equation}\label{Intro-1}
\left\{
\begin{array}{lll}
\partial_t y-\Delta y+a(x,t)y=\chi_{\omega_1}u_1+\chi_{\omega_2}u_2&\mbox{in}&\Omega\times (0,T),\\
y=0&\mbox{on}&\partial\Omega\times (0,T),\\
y(0)=y_0&\mbox{in}&\Omega,
\end{array}\right.
\end{equation}
where $y_0\in L^2(\Omega)$ and $a\in L^\infty(\Omega\times (0,T))$ are
two given functions. For each $u_i\in L^2(0,T;L^2(\Omega)) (i=1,2)$,
(\ref{Intro-1}) has a unique solution in $C([0,T];L^2(\Omega))$, denoted by $y(\cdot;u_1,u_2)$.
Throughout the paper, we will omit the variables $x$ and $t$ for functions of $(x,t)$ and the variable
$x$ for functions of $x$, if there is no risk of causing any confusion.

For each $i=1,2$, we define the following admissible set of controls:
\begin{equation*}
{\mathcal U}_i\triangleq \{u\in L^\infty(0,T;L^2(\Omega)): \|u(t)\|_{L^2(\Omega)}\leq M_i\;\;
\mbox{for a.e.}\;\;t\in (0,T)\},
\end{equation*}
where $M_i$ is a positive constant. Meanwhile, we introduce the following two functionals: For each $i=1,2$,
the functional $J_i: (L^2(0,T;L^2(\Omega)))^2\rightarrow [0,+\infty)$ is defined by
\begin{equation}\label{Intro-2}
J_i(u_1,u_2)\triangleq \|y(T;u_1,u_2)-y_i\|_{L^2(\Omega)},\;\;\mbox{for all}\;\;(u_1,u_2)
\in (L^2(0,T;L^2(\Omega)))^2,
\end{equation}
where $y_i\in L^2(\Omega)$ is given and $y_1\not=y_2$. The problem we consider in this paper is:\\

{\bf(P)}\;\;Does there exist $(u_1^*,u_2^*)\in {\mathcal U}_1\times {\mathcal U}_2$ so that
\begin{equation}\label{20180208-1}
J_1(u_1^*,u_2^*)\leq J_1(u_1,u_2^*)\;\;\mbox{for all}\;\;u_1\in {\mathcal U}_1
\end{equation}
and
\begin{equation}\label{20180208-2}
J_2(u_1^*,u_2^*)\leq J_2(u_1^*,u_2)\;\;\mbox{for all}\;\;u_2\in {\mathcal U}_2?
\end{equation}
We call the problem {\bf(P)} as a two-person nonzero-sum game problem. If the answer to the
problem {\bf(P)} is yes, we call $(u_1^*,u_2^*)$ a Nash equilibrium (or an optimal
strategy pair, or an optimal control pair) of {\bf(P)}. We can understand the problem {\bf(P)} in the following manner:
There are two players executing their strategies and hoping to achieve their goals $y_1$ and $y_2$, respectively.
If the first player chooses the strategy $u_1^*$, then the second player can execute the strategy $u_2^*$
so that $y(T;u_1^*,u_2^*)$ is closer to $y_2$; Conversely, if the second player chooses the
strategy $u_2^*$, then the first player can execute the strategy $u_1^*$ so that $y(T;u_1^*,u_2^*)$ is closer
to $y_1$. Roughly speaking, if one player is deviating from $(u_1^*,u_2^*)$, then the cost functional of this
player would get larger; and there is no information given if both players are deviating from
the Nash equilibrium $(u_1^*,u_2^*)$.

The first main result of this paper is about the existence of a Nash equilibrium of the problem {\bf(P)}.
\begin{Theorem}\label{Intro-3}
The problem {\bf(P)} admits a Nash equilibrium, i.e., there exists
a pair of $(u_1^*,u_2^*)\in {\mathcal{U}}_1\times {\mathcal{U}_2}$
so that (\ref{20180208-1}) and (\ref{20180208-2}) hold.
\end{Theorem}

Differential games were introduced originally by Isaacs (see \cite{Isaacs-1} and \cite{Isaacs-2}).
Since then, lots of researchers were attracted to establish and improve the related theory. Meanwhile,
the theory was applied to a large number of fields. For a comprehensive survey on the differential
game theory, we refer to \cite{Friedman}, \cite{Elliot}, \cite{Basar}, \cite{Evans-Souganidis},
\cite{Yong} and the references therein. It is worthy to mention that, in the vase literature on
game theory, the Kakutani type fixed point theorems were often used to obtain the existence
of Nash equilibria (see, for instance, \cite{Nash}). However, to the best of our knowledge, it
seems the first time to apply this approach to the game problems for heat equation.\\

The second main result of this paper is concerned with the bang-bang property of the Nash
equilibria of the problem {\bf(P)}.
\begin{Theorem}\label{Intro-4}
Let $(u_1^*,u_2^*)\in {\mathcal{U}}_1\times {\mathcal{U}_2}$ be a Nash equilibrium of the problem
{\bf(P)}. Then
\begin{equation*}
\|u_1^*(t)\|_{L^2(\Omega)}=M_1\;\;\mbox{for a.e.}\;\;t\in (0,T)
\end{equation*}
or
\begin{equation*}
\|u_2^*(t)\|_{L^2(\Omega)}=M_2\;\;\mbox{for a.e.}\;\;t\in (0,T).
\end{equation*}
\end{Theorem}

To the best of our knowledge, the studies on bang-bang property are mainly about time optimal control problems. Bang-bang property
is indeed a property of time optimal controls. It is not only important from perspective of applications,
but also very interesting from perspective of mathematics. In some cases, from the bang-bang property of time
optimal controls, one can easily get the uniqueness of time optimal control to this problem (see \cite{Fattorini}).
The bang-bang property may help us to do better numerical analyses and algorithm on time optimal controls in
some cases (see, for instance, \cite{Kunisch-1}, \cite{Kunisch-2} and \cite{Lu}). One of the usual methods to
derive the bang-bang property is the use of the controllability from measurable sets in time (see, for instance, \cite{Wang-2},
\cite{Phung-2}, \cite{Phung-3} and \cite{Apraiz}).
Another usual method to derive the bang-bang property is the use of the Pontryagin Maximum Principle,
together with some unique continuation of the adjoint equation (see, for instance, \cite{Chen} and \cite{Kunisch-3}).

In the next two sections of this paper, we will give the proofs of Theorem~\ref{Intro-3} and
Theorem~\ref{Intro-4}, respectively.

\section{Existence of Nash equilibrium}

In this section, we will prove Theorem~\ref{Intro-3}. Its proof needs the next Kakutani Fixed Point Theorem
quoted from \cite{Aliprantis}.

\begin{Lemma}\label{Intro-5}  Let $S$ be a nonempty, compact and convex subset of
a locally convex Hausdorff space $X$. Let $\Phi: S\mapsto 2^S$ (where $2^S$ denotes the set consisting of
all subsets of $S$)
be a set-valued function
satisfying:
\begin{itemize}
\item[(i)] For each $s\in S$, $\Phi(s)$ is a nonempty and convex subset;
\item[(ii)] Graph$\Phi\triangleq\{(s,z): s\in S\;\mbox{and}\;z\in \Phi(s)\}$ is closed.
\end{itemize}
Then the set of fixed points of $\Phi$ is nonempty and compact, where $s^*\in S$ is called to
be a fixed point of $\Phi$ if $s^*\in \Phi(s^*)$.
\end{Lemma}

\noindent\textbf{\it{Proof of Theorem~\ref{Intro-3}.}}
We first introduce  three set-valued functions
$\Phi_1: {\mathcal{U}_1}\mapsto 2^{{\mathcal{U}_2}}, \Phi_2: {\mathcal{U}_2}\mapsto 2^{{\mathcal{U}_1}}$
and $\Phi: {\mathcal{U}_1}\times {\mathcal{U}_2}\mapsto 2^{{\mathcal{U}_1}\times {\mathcal{U}_2}}$
as follows:
\begin{equation}\label{result-1}
\Phi_1 u_1\triangleq \{u_2\in {\mathcal U}_2: J_2(u_1,u_2)\leq J_2(u_1,v_2)\;\;
\mbox{for all}\;\;v_2\in {\mathcal U}_2\},\;\;u_1\in {\mathcal U}_1,
\end{equation}
\begin{equation}\label{result-2}
\Phi_2 u_2\triangleq \{u_1\in {\mathcal U}_1: J_1(u_1,u_2)\leq J_1(v_1,u_2)\;\;
\mbox{for all}\;\;v_1\in {\mathcal U}_1\},\;\;u_2\in {\mathcal U}_2,
\end{equation}
and
\begin{equation}\label{result-3}
\Phi(u_1,u_2)\triangleq \{({\widetilde u}_1,{\widetilde u}_2): {\widetilde u}_1\in \Phi_2 u_2\;\;\mbox{and}\;\;
{\widetilde u}_2\in \Phi_1 u_1\},\;\;(u_1,u_2)\in {\mathcal U}_1\times {\mathcal U}_2.
\end{equation}
Then we set
\begin{equation*}
X\triangleq (L_w^\infty(0,T;L^2(\Omega)))^2\;\;\mbox{and}\;\;S\triangleq {\mathcal U}_1\times {\mathcal U}_2.
\end{equation*}
It is clear that $X$ is a locally convex Hausdorff space.
The rest of the proof will be carried out by the following four steps.\\

Step 1. We show that $S$ is a nonempty, compact and convex subset of $X$.

This fact can be easily checked. We omit the proofs here.\\

Step 2. We prove that $\Phi(u_1,u_2)$ is nonempty for each $(u_1,u_2)\in S$.

We arbitrarily fix $(u_1,u_2)\in S$. According to (\ref{result-1})-(\ref{result-3}), it suffices to
show that $\Phi_1 u_1$ and $\Phi_2 u_2$ are nonempty. For this purpose, we introduce the following
auxiliary optimal control problem:
\begin{equation*}
{\rm\bf (P_{au})}\;\;\;\;\;\;\;\inf_{v_2\in {\mathcal U}_2} J_2(u_1,v_2).
\end{equation*}
Let
\begin{equation}\label{result-4}
d\triangleq\inf_{v_2\in {\mathcal U}_2} J_2(u_1,v_2).
\end{equation}
It is obvious that $d\geq 0$. Let $\{v_{2,n}\}_{n\geq 1}\subseteq {\mathcal U}_2$ be a minimizing sequence
so that
\begin{equation}\label{result-5}
d=\lim_{n\rightarrow \infty} J_2(u_1,v_{2,n}).
\end{equation}
On one hand, since $\|v_{2,n}\|_{L^\infty(0,T;L^2(\Omega))}\leq M_2$,
there exists a subsequence of $\{n\}_{n\geq 1}$, still denoted by itself, and
$v_{2,0}\in {\mathcal U}_2$, so that
\begin{equation}\label{result-5:1}
v_{2,n}\rightarrow v_{2,0}\;\;\mbox{weakly star in}\;\;L^\infty(0,T;L^2(\Omega)).
\end{equation}
On the other hand, we denote that $z_n(\cdot)\triangleq y(\cdot;u_1,v_{2,n})-y(\cdot;u_1,v_{2,0})$.
According to (\ref{Intro-1}), it is clear that
\begin{equation}\label{result-5:2}
\left\{
\begin{array}{lll}
\partial_t z_n-\Delta z_n+a(x,t)z_n=\chi_{\omega_2}(v_{2,n}-v_{2,0})&\mbox{in}&\Omega\times (0,T),\\
z_n=0&\mbox{on}&\partial\Omega\times (0,T),\\
z_n(0)=0&\mbox{in}&\Omega.
\end{array}\right.
\end{equation}
By $L^2$-theory for parabolic equation (see \cite{Evans}), we obtain that
\begin{equation}\label{result-5:3}
\|z_n\|_{L^2(0,T;H^2(\Omega)\cap H_0^1(\Omega))}+\|\partial_t z_n\|_{L^2(0,T;L^2(\Omega))}
\leq C\|v_{2,n}-v_{2,0}\|_{L^2(0,T;L^2(\Omega))},
\end{equation}
where $C>0$ is a constant independent of $n$. It follows from
(\ref{result-5:1}) and (\ref{result-5:3}) that there exists a subsequence of
$\{n\}_{n\geq 1}$, still denoted by itself, and $z\in L^2(0,T;H^2(\Omega)\cap H_0^1(\Omega))
\cap W^{1,2}(0,T;L^2(\Omega))$, so that
\begin{equation}\label{result-5:4}
\begin{array}{ll}
z_n\rightarrow z&\mbox{weakly in}\;\; L^2(0,T;H^2(\Omega)\cap H_0^1(\Omega))
\cap W^{1,2}(0,T;L^2(\Omega))\\
&\mbox{and strongly in}\;\;C([0,T];L^2(\Omega)).
\end{array}
\end{equation}
Passing to the limit for $n\rightarrow \infty$ in (\ref{result-5:2}), by
(\ref{result-5:1}) and (\ref{result-5:4}), we obtain that $z=0$. Hence,
\begin{equation}\label{result-6}
y(T;u_1,v_{2,n})\rightarrow y(T;u_1,v_{2,0})\;\;\mbox{strongly in}\;\;L^2(\Omega).
\end{equation}
It follows from (\ref{result-5}), (\ref{Intro-2}) and (\ref{result-6}) that
\begin{equation}\label{result-7}
d=J_2(u_1,v_{2,0}).
\end{equation}
Noting that $v_{2,0}\in {\mathcal U}_2$, by (\ref{result-4}), (\ref{result-7}) and (\ref{result-1}),
we obtain that $v_{2,0}\in \Phi_1 u_1$. This implies that $\Phi_1 u_1\not=\emptyset$. In the same way,
we also have that $\Phi_2 u_2\not=\emptyset$.\\

Step 3. We show that $\Phi(u_1,u_2)$ is a convex subset of ${\mathcal U}_1\times {\mathcal U}_2$ for
each $(u_1,u_2)\in {\mathcal U}_1\times {\mathcal U}_2$.

We arbitrarily fix $(u_1,u_2)\in {\mathcal U}_1\times {\mathcal U}_2$.
According to (\ref{result-1})-(\ref{result-3}), it suffices to prove that
$\Phi_1 u_1$ is a  convex subset of ${\mathcal U}_2$. The convexity of $\Phi_2 u_2$ can
be similarly proved. For this purpose, we arbitrarily fix ${\widetilde u}_2, {\widehat u}_2\in \Phi_1 u_1$.
By (\ref{result-1}),
we get that
\begin{equation}\label{result-8}
{\widetilde u}_2,\;\;{\widehat u}_2\in {\mathcal U}_2,
\end{equation}
\begin{equation}\label{result-9}
J_2(u_1,{\widetilde u}_2)\leq J_2(u_1,v_2)\;\;\mbox{and}\;\;
J_2(u_1,{\widehat u}_2)\leq J_2(u_1,v_2)\;\;\mbox{for each}\;\;v_2\in {\mathcal U}_2.
\end{equation}
For any $\lambda\in [0,1]$, by (\ref{Intro-2}) and (\ref{Intro-1}), we have that
\begin{eqnarray*}
J_2(u_1,\lambda {\widetilde u}_2+(1-\lambda){\widehat u}_2)&=&
\|y(T;u_1,\lambda {\widetilde u}_2+(1-\lambda){\widehat u}_2)-y_2\|_{L^2(\Omega)}\\
&=&\|\lambda[y(T;u_1,{\widetilde u}_2)-y_2]+(1-\lambda)[y(T;u_1,{\widehat u}_2)-y_2]\|_{L^2(\Omega)}\\
&\leq&\lambda J_2(u_1,{\widetilde u}_2)+(1-\lambda)J_2(u_1,{\widehat u}_2).
\end{eqnarray*}
This, along with (\ref{result-8}) and (\ref{result-9}), yields that
\begin{equation*}
\lambda {\widetilde u}_2+(1-\lambda){\widehat u}_2\in {\mathcal U}_2
\end{equation*}
and
\begin{equation*}
J_2(u_1,\lambda{\widetilde u}_2+(1-\lambda){\widehat u}_2)\leq J_2(u_1,v_2)\;\;\mbox{for each}\;\;v_2\in {\mathcal U}_2,
\end{equation*}
which indicate that $\lambda{\widetilde u}_2+(1-\lambda){\widehat u}_2\in \Phi_1 u_1$ (see (\ref{result-1})). Hence, $\Phi_1 u_1$
is a convex subset of ${\mathcal U}_2$. \\

Step 4. We prove that Graph$\Phi$ is closed.

It suffices to show that if $(u_{n,1},u_{n,2})\in {\mathcal U}_1\times {\mathcal U}_2,
{\widetilde u}_{n,1}\in \Phi_2 u_{n,2}$, ${\widetilde u}_{n,2}\in \Phi_1 u_{n,1}$,
$(u_{n,1},u_{n,2})\rightarrow (u_1,u_2)$ in $X$
and $({\widetilde u}_{n,1},{\widetilde u}_{n,2})\rightarrow ({\widetilde u}_1,{\widetilde u}_2)$  in
$X$, then
\begin{equation}\label{20180208-3}
(u_1,u_2)\in {\mathcal U}_1\times {\mathcal U}_2,\; {\widetilde u}_1\in \Phi_2 u_2\;\;
\mbox{and}\;\;{\widetilde u}_2\in \Phi_1 u_1.
\end{equation}
Indeed, on one hand, by (\ref{result-1}) and (\ref{result-2}), we can easily check that
\begin{equation}\label{result-10}
(u_1,u_2)\in {\mathcal U}_1\times {\mathcal U}_2,\; {\widetilde u}_1\in {\mathcal U}_1\;\;
\mbox{and}\;\;{\widetilde u}_2\in {\mathcal U}_2.
\end{equation}
On the other hand, according to ${\widetilde u}_{n,1}\in \Phi_2 u_{n,2}$, (\ref{result-2}) and (\ref{Intro-2}), it is obvious that
\begin{equation}\label{result-11}
\|y(T;{\widetilde u}_{n,1},u_{n,2})-y_1\|_{L^2(\Omega)}\leq \|y(T;v_1,u_{n,2})-y_1\|_{L^2(\Omega)}\;\;\mbox{for each}\;\;
v_1\in {\mathcal U}_1.
\end{equation}
Since $({\widetilde u}_{n,1},u_{n,2})\rightarrow ({\widetilde u}_1,u_2)$ weakly star in $(L^\infty(0,T;L^2(\Omega)))^2$,
by similar arguments as those to get (\ref{result-6}), there
exists a subsequence of $\{n\}_{n\geq 1}$, still denoted by itself, so that
\begin{equation}\label{result-12}
(y(T;{\widetilde u}_{n,1},u_{n,2}),y(T;v_1,u_{n,2}))\rightarrow (y(T;{\widetilde u}_1,u_2),y(T;v_1,u_2))\;\;\mbox{strongly in}\;\;(L^2(\Omega))^2.
\end{equation}
Passing to the limit for $n\rightarrow \infty$ in (\ref{result-11}), by (\ref{result-12}),
we get that
\begin{equation*}
\|y(T;{\widetilde u}_1,u_2)-y_1\|_{L^2(\Omega)}\leq \|y(T;v_1,u_2)-y_1\|_{L^2(\Omega)}\;\;\mbox{for each}\;\;v_1\in {\mathcal U}_1.
\end{equation*}
This, together with (\ref{Intro-2}), (\ref{result-2}) and the second conclusion in (\ref{result-10}), implies that
${\widetilde u}_1\in \Phi_2 u_2$. Similarly, ${\widetilde u}_2\in \Phi_1 u_1$. Hence, (\ref{20180208-3}) follows.\\

Step 5. We finish the proof.

According to Steps 1-4 and Lemma~\ref{Intro-5}, there exists a pair of $(u_1^*,u_2^*)\in {\mathcal U}_1\times {\mathcal U}_2$
so that $(u_1^*,u_2^*)\in \Phi(u_1^*,u_2^*)$,
which, combined with (\ref{result-1})-(\ref{result-3}), indicates that
$(u_1^*,u_2^*)$ is a Nash equilibrium of the problem {\bf(P)}.\\

In summary, we end the proof of Theorem~\ref{Intro-3}.\hspace{66mm}$\Box$

\section{Bang-bang property}
\noindent\textbf{\it{Proof of Theorem~\ref{Intro-4}.}} Let $(u_1^*,u_2^*)\in {\mathcal U}_1\times {\mathcal U}_2$ be a Nash equilibrium of the problem {\bf(P)},
i.e., (\ref{20180208-1}) and (\ref{20180208-2}) hold.
We arbitrarily fix $u_1\in {\mathcal U}_1$ and $\lambda\in (0,1)$. Set
$u_{1,\lambda}\triangleq u_1^*+\lambda (u_1-u_1^*)$. It is obvious that $u_{1,\lambda}\in {\mathcal U}_1$.
Then by (\ref{20180208-1}), we get that
\begin{equation*}
J_1(u_1^*,u_2^*)\leq J_1(u_{1,\lambda},u_2^*),
\end{equation*}
i.e., $\|y(T;u_1^*,u_2^*)-y_1\|_{L^2(\Omega)}\leq \|y(T;u_{1,\lambda},u_2^*)-y_1\|_{L^2(\Omega)}$ (see (\ref{Intro-2})).
From the latter it follows that
\begin{equation}\label{result-14}
\langle y(T;u_1^*,u_2^*)-y_1,z(T)\rangle_{L^2(\Omega)}\geq 0,
\end{equation}
where $z\in C([0,T];L^2(\Omega))$ is the unique solution to the equation
\begin{equation}\label{result-15}
\left\{
\begin{array}{lll}
\partial_t z-\Delta z+a(x,t)z=\chi_{\omega_1} (u_1-u_1^*)&\mbox{in}&\Omega\times (0,T),\\
z=0&\mbox{on}&\partial\Omega\times (0,T),\\
z(0)=0&\mbox{in}&\Omega.
\end{array}\right.
\end{equation}
Let $\varphi\in C([0,T];L^2(\Omega))$ be the unique solution to the equation
\begin{equation}\label{result-16}
\left\{
\begin{array}{lll}
\partial_t \varphi+\Delta \varphi-a(x,t)\varphi=0&\mbox{in}&\Omega\times (0,T),\\
\varphi=0&\mbox{on}&\partial\Omega\times (0,T),\\
\varphi(T)=y_1-y(T;u_1^*,u_2^*)&\mbox{in}&\Omega.
\end{array}\right.
\end{equation}
Multiplying (\ref{result-15}) by $\varphi$ and integrating it over $\Omega\times (0,T)$,
after some calculations, by (\ref{result-16}), we obtain that
\begin{equation*}
\int_0^T \langle \chi_{\omega_1}\varphi(t),u_1(t)-u_1^*(t)\rangle_{L^2(\Omega)}\mathrm dt
=\langle z(T),y_1-y(T;u_1^*,u_2^*)\rangle_{L^2(\Omega)}.
\end{equation*}
This, along with (\ref{result-14}), implies that
\begin{equation}\label{result-17}
\int_0^T \langle \chi_{\omega_1}\varphi(t),u_1(t)-u_1^*(t)\rangle_{L^2(\Omega)}\mathrm dt\leq 0\;\;\mbox{for each}\;\;
u_1\in {\mathcal U}_1.
\end{equation}
Since $L^2(\Omega)$ is separable, there exists a countable subset ${\mathbb U}_0=\{v_\ell\}_{\ell\geq 1}$
so that ${\mathbb U}_0$ is dense in $\mathbb U=\{u\in L^2(\Omega): \|u\|_{L^2(\Omega)}\leq M_1\}$.
For each $v_\ell\in {\mathbb U}_0$, we define the function
\begin{equation*}
g_\ell(t)\triangleq \langle \chi_{\omega_1}\varphi(t),v_\ell-u_1^*(t)\rangle_{L^2(\Omega)},\;\;t\in (0,T).
\end{equation*}
Then $g_\ell(\cdot)\in L^1(0,T)$. Thus, there exists a measurable set $E_\ell\subseteq (0,T)$ with $|E_\ell|=T$,
so that any point in $E_\ell$ is a Lebesgue point of $g_\ell(\cdot)$. Namely,
\begin{equation*}
\lim_{\delta\rightarrow 0^+}\frac{1}{\delta}\int_{t-\delta}^{t+\delta}|g_\ell(s)-g_\ell(t)|\mathrm ds=0\;\;
\mbox{for each}\;\;t\in E_\ell.
\end{equation*}
Now, for any $t\in E_\ell$ and $\delta>0$, we define
\begin{equation*}
u_\delta(s)\triangleq\left\{
\begin{array}{lll}
u_1^*(s)&\mbox{if}&s\in (0,T)\setminus (t-\delta,t+\delta),\\
v_\ell&\mbox{if}&s\in (0,T)\cap (t-\delta,t+\delta).
\end{array}\right.
\end{equation*}
Then, by (\ref{result-17}), we get that
\begin{equation*}
\int_{(0,T)\cap (t-\delta,t+\delta)}\langle \chi_{\omega_1}\varphi(s),v_\ell-u_1^*(s)\rangle_{L^2(\Omega)}\mathrm ds\leq 0.
\end{equation*}
Dividing the above inequality by $\delta>0$ and then sending $\delta\rightarrow 0$, we obtain that
$g_\ell(t)\leq 0$. From this, we see that for all $t\in \bigcap_{\ell\geq 1} E_\ell$ and
$v_\ell\in {\mathbb U}_0$,
\begin{equation}\label{result-18}
\langle \chi_{\omega_1}\varphi(t),v_\ell\rangle_{L^2(\Omega)}\leq
\langle \chi_{\omega_1}\varphi(t),u_1^*(t)\rangle_{L^2(\Omega)}.
\end{equation}
Since ${\mathbb U}_0$ is countable and dense in ${\mathbb U}$, by (\ref{result-18}), we have that
$|\bigcap_{\ell\geq 1} E_\ell|=T$ and that for a.e. $t\in (0,T)$,
\begin{equation*}
\langle \chi_{\omega_1} \varphi(t),u\rangle_{L^2(\Omega)}\leq
\langle \chi_{\omega_1}\varphi(t),u_1^*(t)\rangle_{L^2(\Omega)}
\;\;\mbox{for all}\;\;u\in {\mathbb U}.
\end{equation*}
From these we obtain that
\begin{equation}\label{result-19}
\max_{\|u\|_{L^2(\Omega)}\leq M_1}\langle \chi_{\omega_1}\varphi(t),u\rangle_{L^2(\Omega)}=
\langle \chi_{\omega_1}\varphi(t),u_1^*(t)\rangle_{L^2(\Omega)}\;\;\mbox{for a.e.}\;\;t\in (0,T).
\end{equation}
Similarly, we have that
\begin{equation}\label{result-20}
\max_{\|u\|_{L^2(\Omega)}\leq M_2}\langle \chi_{\omega_2}\psi(t),u\rangle_{L^2(\Omega)}=
\langle \chi_{\omega_2}\psi(t),u_2^*(t)\rangle_{L^2(\Omega)}\;\;\mbox{for a.e.}\;\;t\in (0,T),
\end{equation}
where $\psi\in C([0,T];L^2(\Omega))$ is the unique solution to the equation
\begin{equation}\label{result-21}
\left\{
\begin{array}{lll}
\partial_t \psi+\Delta \psi-a(x,t)\psi=0&\mbox{in}&\Omega\times (0,T),\\
\psi=0&\mbox{on}&\partial\Omega\times (0,T),\\
\psi(T)=y_2-y(T;u_1^*,u_2^*)&\mbox{in}&\Omega.
\end{array}\right.
\end{equation}
Noting that $y_1\not=y_2$, by (\ref{result-16}), (\ref{result-21}),
the unique continuation estimate at one time (see \cite{Phung}), and the backward uniqueness of the linear
 heat equation (see \cite{Bardos}), we obtain that
\begin{equation*}
\chi_{\omega_1}\varphi(t)\not=0\;\;\mbox{a.e.}\;\;t\in (0,T)
\;\;\mbox{or}\;\;\chi_{\omega_2}\psi(t)\not=0\;\;\mbox{a.e.}\;\;t\in (0,T).
\end{equation*}
These, together with (\ref{result-19}) and (\ref{result-20}), imply that
\begin{equation*}
u_1^*(t)=M_1\frac{\chi_{\omega_1}\varphi(t)}{\|\chi_{\omega_1}\varphi(t)\|_{L^2(\Omega)}}
\;\;\mbox{a.e.}\;\;t\in (0,T)
\end{equation*}
or
\begin{equation*}
u_2^*(t)=M_2\frac{\chi_{\omega_2}\psi(t)}{\|\chi_{\omega_2}\psi(t)\|_{L^2(\Omega)}}
\;\;\mbox{a.e.}\;\;t\in (0,T).
\end{equation*}
Hence,
\begin{equation*}
\|u_1^*(t)\|_{L^2(\Omega)}=M_1\;\;\mbox{a.e.}\;\;t\in (0,T)\;\;
\mbox{or}\;\;\|u_2^*(t)\|_{L^2(\Omega)}=M_2\;\;\mbox{a.e.}\;\;t\in (0,T).
\end{equation*}

In summary, we finish the proof of Theorem~\ref{Intro-4}.\hspace{63mm}$\Box$




\end{document}